\documentclass[12pt]{article}
\usepackage{amsmath}
\usepackage{amssymb}
\usepackage{amsmath,amssymb,amsthm,amsfonts}
{\theoremstyle{definition}

}

\newcommand{\ph}{\phantom{a}}
\newcommand{\phh}{\phantom{aaa}}

\newcommand{\no}{\eqno}

\newcommand{\al}{\alpha}

\newcommand{\et}{\eta}

\newcommand{\la}{\lambda}
\newcommand{\La}{\Lambda}
\newcommand{\fr}{\frac}
\newcommand{\sist}[2]{\left\{
\begin{array}{l}
{#1}\\
\ph\\
{#2}
\end{array}
\right.}

\begin{document}

MSC 68W05, 68W99, 12G99

\vskip 20pt

\centerline{\bf On the localization of roots of polynomials}

\vskip 10 pt

\centerline{\bf G. A. Grigorian}

\vskip 10 pt

\centerline{0019 Armenia c. Yerevan, str. M. Bagramian 24/5}
\centerline{Institute of Mathematics of NAS of Armenia}
\centerline{E - mail: mathphys2@instmath.sci.am, \ph phone: 098 62 03 05, \ph 010 35 48 61}

\vskip 20 pt

\noindent
Abstract. In this article we use a method of finding the index of a complex-valued function by determined number of arithmetic operations to describe an algorithm of localization of roots of square-free polynomials. We give an estimation of the number of arithmetic operations for the described algorithm.

\vskip 20 pt


\vskip 20 pt

{\bf  1. Introduction}. The main theorem of Algebra states that every polynomial of degree $n$
$$
P(z)\equiv a_0 z^n +\dots + a_n \phh (a_j\in \mathbb{C}, \ph j=\overline{0,n}, \ph a_0\ne 0)
$$
is representable in the form
$$
P(z) = a_0(z - z_1)^{n_1}\dots (z - z_d)^{n_d}, \ph n_j\in \mathbb{N}, \ph j=\overline{1,d}, \ph n_1 + \dots + n_d= n.
$$
The numbers $z_1, \dots, z_d$ are called the roots (the zeroes) of $P(z)$. It is well  known that for $n> 4$ the roots of $P(z)$ are not expressed by its coefficients via radicals. Therefore arises the question how to find all roots with their multiplicity of $P(z)$  numerically with given exactness. More exactly, there arises the question: how for every enough small $\varepsilon > 0$  to find the discs of radius $\varepsilon$ in $\mathbb{C}$, each of which contain exactly one root of $P(z)$. This problem is not solved yet and   many works are devoted to it (see [1-4] and cited works therein).

In this paper we give an algorithm to find the mentioned above discs of radius $\varepsilon$ for every $\varepsilon > 0$ by determined number (depending on $\varepsilon$) of arithmetic operations. We give an estimation of  the  number of these arithmetic operations via $\varepsilon$ in the case of square-free $P(z)$ (a necessary and sufficient condition for $P(z)$ to be square-free is that $Resultant (P,P) \ne 0$, see below).

\vskip 10pt

{\bf 2. Auxiliary propositions.} Along with $P(z)$ consider the polynomial
$$
Q(z)\equiv b_0 z^m +\dots + b_m \phh (a_j\in \mathbb{C}, \ph j=\overline{0,m}, \ph a_0\ne 0)
$$
Important tools for studying the algebra of polynomials  is the resultant of polynomials, for example, of $P(z)$ and $Q(t)$ (see [5, p.126])
$$
Resultant (P,Q) \equiv \begin{vmatrix} a_0&a_1&\dots a_n&\phantom{a}&\phantom{a}\\ \phantom{a}&a_0&a_1&\dots&a_n& .\\
\dots&\dots&\dots &\dots&\dots& .\\
\phantom{a}& \phantom{a}&a_0&a_1&\dots&.& a_n\\
b_0&b_1&\dots b_m&\phantom{a}&\phantom{a}\\ \phantom{a}&b_0&b_1&\dots&b_m& .\\
\dots&\dots&\dots &\dots&\dots& .\\
\phantom{a}& \phantom{a}&b_0&b_1&\dots&.& b_m
 \end{vmatrix}
$$
and discriminant of a Polynomial, say $P(z)$:
$$
\mathcal{D}(P) \equiv a_0^{2n-2}\prod \limits_{j <k}(\eta_j - \eta_k),
$$
where $\eta_j, \ph j=\overline{1,n}$ are roots of $P(z)$, counted as many times as their multiplicity.
It is well known that (see [5, p. 130])
$$
Resultant (P,P') = \pm a_0 \mathcal{D}(P)
$$
It follows from here that $P(z)$ is square-free if and only if $Resultant(P,P') \ne 0$. This condition one cam check effectively by finite number of arithmetic operations with the coefficients of $P(z)$ and $P'(z)$.

\vskip 10pt

{\bf Theorem 2.1 ([1, Theorem 2, the Cauchy inequality]).} {\it Any root $z$ of $P$ satisfies
$$
|z| < \fr{\max\{|a_1|,\dots,|a_n|\}}{|a_0|} \stackrel{def}{=} \mathcal{R}_0.
$$
}

\phantom{aaaaaaaaaaaaaaaaaaaaaaaaaaaaaaaaaaaaaaaaaaaaaaaaaaaaaaaaaaaaaaaaaaaaaaa}$\Box$

\vskip 10pt

{\bf Corollary 2.1 ([1, Corollary]).}  {\it Let $a_n \ne 0$. Then any root $z$ of $P$ satisfies
$$
|z| > \fr{|a_n|}{|a_n| + \max\{|a_0|, \dots, |a_{n-1}|\}} \stackrel{def}{=} \rho_0.
$$
}

\phantom{aaaaaaaaaaaaaaaaaaaaaaaaaaaaaaaaaaaaaaaaaaaaaaaaaaaaaaaaaaaaaaaaaaaaaaa}$\Box$

\vskip 10pt

We set $||P|| \equiv \sqrt{|a_0|^2+|a_1|^2+\dots+|a_n|^2}, \ph sep(P) \equiv \min\limits_{j\ne k}|z_j - z_k|$, where $z_j, \ph j=\overline{1,n}$ are all roots of $P(z)$.

\vskip 10pt

{\bf Theorem 2.2 ([1, Theorem 5]).} {\it Let $P(z)$ be square-free. Then
$$
sep(P) > \sqrt{3}n^{-\fr{n+2}{2}}|\mathcal{D}(P)|^{1/2}||P||^{1-n}\stackrel{def}{=}\varepsilon_0.
$$
}

\phantom{aaaaaaaaaaaaaaaaaaaaaaaaaaaaaaaaaaaaaaaaaaaaaaaaaaaaaaaaaaaaaaaaaaaaaaa}$\Box$

Denote by $D(A,r)$ a disc of a radius $r$  with the center in the point $A \in \mathbb{C}$:
$$
D(A,r)\equiv \{z \in \mathbb{C} \ph : \ph  |z - A| \le r\},
$$
and denote by $C_r$ a circle of radius $r$ with the center in $A$:
$$
C_r\equiv \{z \in \mathbb{C} \ph : \ph  |z| = r\}.
$$
Finally, denote by $K(a,b)$ ($0< a < b$) an annulus of radiuses $a$ and $b$ with the center in $0$, that is
$$
K(a,b) \equiv \{z \in \mathbb{C} \ph : \ph a \le |z| \le b\}.
$$
and denote by $\stackrel{_\circ}{K}(a,b)$ the inner part of $K(a,b)$. The number $b-a$ we will call the width of $K(a,b)$ (of $\stackrel{_\circ}{K}(a,b)$).

{\bf Lemma 2.1.} {\it Let $K(a,b)$ be an annulus free of roots of $P$ and $P'$. Then
$$
\min\limits_{|z| = a + k \fr{b-a}{3}}\left|\fr{P'(z)|}{P(z)}\right| \ge \frac{n}{(2\mathcal{R}_0)^n}\Bigl(\fr{b-a}{3}\Bigr)^{n-1}, \phh k=1,2. \eqno (2.1)
$$
}

Proof. Let $P'(z) = a_0 n (z - \xi_1)^{l_1}\dots (z - \xi_s)^{l_s}, \ph l_j \in \mathbb{N}, \ph \ph j=\overline{1,s}, \ph l_1 + \dots + l_s = n-1.$ Then
$$
\left|\fr{P'(z)|}{P(z)}\right| = \frac{n |z - \xi_1|^{l_1}\dots |z - \xi_s|^{l_s}}{|z - z_1|\dots|z -z_n|}. \eqno (2.2)
$$
Since $K(a,b)$ is free of roots of $P(z)$ and $P'(z)$ we have $|z - \xi_j| \ge \fr{b-a}{3}, \ph |z| = a + k\fr{b-a}{3}, \ph k=1,2, \ph j=\overline{1,s},$ and by Theorem 2.1  $|z - z_j| \le 2 \mathcal{R}_0, \ph |z| = a + k\fr{b-a}{3}, \ph k=1,2, \ph j=\overline{1,n}.$ From here and from (2.2) it follows (2.1). The lemma is proved.

Consider the function
$$
f(r,\theta) \equiv \fr{P'(r e^{i\theta})}{P(r e^{i\theta})}, \ph r \ge 0, \ph \theta \in \mathbb{R}, \ph  r e^{i\theta} \ne z_j, \ph j=\overline{1.n}.
$$

{\bf Lemma 2.2.} {\it Let $K(a,b)$ be an annulus free of roots of $P(z)$. Then
$$
\left | \fr{\partial} {\partial \theta}f(r,\theta) \right| \le \frac{9 n b}{(b -a)^2}, \ph r  = a + k\fr{b-a}{3}, \ph k=1,2, \ph \theta \in \mathbb{R}. \eqno (2.3)
$$
}

Proof. We have
$$
f(r,\theta) = \fr{1}{r e^{i\theta} - z_1} + \dots + \fr{1}{r e^{i\theta} - z_n}, \ph r \ge 0, \ph \theta \in \mathbb{R}, \ph r e^{i\theta} \ne  z_j, \ph j=\overline{1,n}.
$$
Then
$$
\fr{\partial} {\partial\theta}f(r,\theta) = - i r \left(\fr{1}{(r e^{i\theta} - z_1)^2} + \dots + \fr{1}{(r e^{i\theta} - z_n)^2}\right), \ph r \ge 0, \ph \theta \in \mathbb{R}, \no (2.4)
$$
$r e^{i\theta} \ne  z_j, \ph j=\overline{1,n}.$ Since $K(a,b)$ is free of roots of $P(z)$ we have $|r e^{i\theta} - z_j| \ge \fr{b-a}{3}, \ph j=\overline{1.n}, \ph r = a + k \fr{b-a}{3}, \ph k=1,2.$ This together with (2.4) implies (2.3). The lemma is proved.

Denote by $KS(a,b,\lambda,\mu)$ a segment of an annulus $K(a,b)$, bordered by the angles $\theta =\lambda < \theta = \mu$, that is
$$
K(a,b,\lambda,\mu) \equiv \{z = r e^{i\theta} \ph : \ph a \le r \le b, \ph \lambda \le \theta \le \mu\}.
$$
and denote by $\stackrel{_\circ}{K}S(a,b,\la,\mu)$ the inner part of $KS(a,b,\lambda,\mu)$. The number $b-a$ we will call the width of  $K(a,b,\lambda,\mu)$ (of $\stackrel{_\circ}{K}S(a,b,\la,\mu)$) and the number $\mu-\la$ we will call the length of $K(a,b,\lambda,\mu)$ (of $\stackrel{_\circ}{K}S(a,b,\la,\mu)$).

\vskip 10pt

{\bf Lemma 2.3.} {\it Let $0 < c < a <b, \ph \mu - \lambda < \pi$ and let $KS(a -c, b+c, \lambda,\mu)$ be free of roots of $P(z)$ and $P'(z)$. Then
$$
\min\limits_{a \le r \le b}\left|f\left(r,\fr{\lambda + \mu}{2}\right)\right| \ge \fr{n(\min\{c, a\tan\fr{\mu -\lambda}{2}\})^{n-1}}{(2\mathcal{R}_0)^n}. \no (2.5)
$$
}

Proof. By (2.2) we have
$$
\left| f\left(r,\fr{\la + \mu}{2}\right)\right| = \fr{n|r e^{i\fr{\la +\mu}{2}} - \xi_1|^{l_1}\dots |r e^{i\fr{\la +\mu}{2}} - \xi_s|^{l_s}}{|r e^{i\fr{\la +\mu}{2}} - z_1|\dots |r e^{i\fr{\la +\mu}{2}} - z_n|}, \ph l_1 + \dots + l_s = n-1. \no (2.6)
$$
We set $c_1\equiv \min\{c, \sin\fr{\mu - \la}{2}\}$. Chose $\la_1$ and $\mu_1$ such that $\la \le \la_1 \le \mu_1\le \mu, \ph \mu_1 - \la_1 < \pi$ and $2\tan \fr{\mu_1 - \la_1}{2} =~ c_1\le c$. Then, obviously (see Figure 1),
$$
KS(a-c_1, b+c_1, \la_1,\mu_1) \subset KS(a-b,b+c,\la,\mu)
$$
We set $A_1\equiv a e^{i\fr{\la+\mu}{2}}, \ph A_2\equiv b e^{i\fr{\la+\mu}{2}}$. Let $[B_1B_4]$ and $[B_2B_3]$ be intervals, which are tangent to the discs $D(A_1,c_1)$ and $D(A_2,c_1)$ at the points $B_1$ and  $B_2$  respectively such that $,B_1\equiv a e^{i\la_1}, \ph B_2\equiv a e^{i\mu_1} \in C_a, \ph B_3, B_4 \in C_b, \ph |A_1 - B_1| = |A_1 - B_2| = |A_2- B_3| = |A_2-B_4| = c_1 = \min\{c,2\tan\fr{\mu_1-\la_1}{2}\}.$ Then the obtained curve "parallelogram" $\Pi\equiv B_1D_1B_2B_3D_2B_4B_1$(see Figure 1) is a subset of $KS(a-c_1,b+c_1,\la_1,\mu_1)$. Hence, $\Pi$ is free of roots of $P(z)$ and$P'(z)$.
Obviously, the distance of every $z\in A_1A_2 \equiv\{r e^{\fr{\la+\mu}{2}} \ph : \ph a \le r \le b\}$ to the border of
$\Pi$ is $c_1=\min\{c, a\tan\fr{\mu_1-\la_1}{2}\}.$ Therefore,
$$
|z - \xi_j| \ge c_1, \phh j=\overline{1,s}, \phh z \in A_1A_2. \no (2.7)
$$
\vskip 50pt

\begin{picture}(120,190)
\put(260,35){\line(2,3){110}}
\put(20,200){\line(2,-3){100}}
\put(200,47){\line(1,4){47}}
\put(130,233){\line(1,-4){47}}
\put(150, 0) {$Figure \ph 1.$}
\qbezier[360](20,200)(200,280)(370,200)
\qbezier[360](50,150)(200,220)(340,150)
\qbezier[360](80,100)(190,170)(305,100)
\qbezier[360](110,50)(190,120)(270,50)
\qbezier[160](100,80)(190,136)(290,75)
\put(220, 120){$B_1$}
\put(140, 120){$B_2$}
\put(141, 185){$B_3$}
\put(221, 185){$B_4$}
\put(180, 120){$A_1$}
\put(180, 190){$A_2$}

\put(180, 95){$D_1$}
\put(180, 210){$D_2$}
\put(300, 90){$ae^{i\la}$}
\put(60, 100){$ae^{i\mu}$}
\put(340, 140){$be^{i\la}$}
\put(35, 140){$be^{i\mu}$}
\put(270, 40){$(a-c)e^{i\la}$}
\put(180, 65){$(a-c)e^{i\la_1}$}
\put(112, 82){$(a-c)e^{i\mu_1}$}
\put(70, 40){$(a-c)e^{i\mu}$}

\put(292,72){$(a-c_1)e^{i\la}$}

\put(370, 190){$(b+c)e^{i\la}$}
\put(230, 240){$(b+c)e^{i\la_1}$}
\put(-30, 190){$(b+c)e^{i\mu}$}
\put(110, 240){$(b+c)e^{i\mu_1}$}
\put(187,135){\line(0,1){50}}
\put(156,132){\line(0,1){50}}
\put(220,133){\line(0,1){51}}
\qbezier[360](157,131)(162,110)(180,107)
\qbezier[360](180,107)(202,101)(220,130)

\qbezier[360](157,182)(163,208)(190,207)
\qbezier[360](190,207)(209,206)(220,185)

\end{picture}

\vskip 30pt

By virtue of Theorem 2.1 we have
$$
|z-z_j|\le 2\mathcal{R}_0, \phh z \in A_1A_2, \phh j=\overline{1,n}.
$$
This together with (2.6) and (2.7) implies (2.5). The lemma is proved.

\vskip 10pt

{\bf Lemma 2.4.} {\it Let $0<c<a<b,\phh \mu-\la<\pi$ and let $KS(a-c,b+c,\la,\mu)$ be free of roots of $P(z)$. Then
$$
\left|\fr{\partial f\left(r,\fr{\la+\mu}{2}\right)}{\partial r}\right| \le \fr{n}{(\min\{c,a\tan\fr{\mu-\la}{2}\})^2}, \phh a\le r\le b. \no (2.8)
$$
}

Proof. We have
$$
\left|\fr{\partial f\left(r,\fr{\la+\mu}{2}\right)}{\partial r}\right| \le \fr{1}{|r e^{i\fr{\la+\mu}{2}} - z_1|^2}+\dots+\fr{1}{|r e^{i\fr{\la+\mu}{2}} - z_n|^2}, \phh a\le r \le b. \no (2.9)
$$
It was shown in the proof of Lemma 2.3 that $\Pi\subset KS(a-c,b+c,\la,\mu)$.  Then since $KS(a-c,b+c,\la,\mu)$ is free of roots of $P(z)$, the "curve parallelogram" $\Pi$ is also free of roots of $P(z)$. Then
$$
|re^{i\fr{\la+\mu}{2}} - z_j|\ge \min\Bigl\{c, a \tan\fr{\mu-\la}{2}\Bigr\}, \phh j=\overline{1,n}, \phh  a\le r \le b.
$$
This together with (2.9) implies (2.8). The lemma is proved.

We set $\theta_{j+1}\equiv \theta_0 + j \fr{2\pi}{N}, \ph j=\overline{1,N-1}, \ph r_j\equiv j\fr{\mathcal{R}_0}{N}, \ph j=\overline{0,N}$. By the Lagrange formula we have
$$
f(r,\theta) - f(r,\theta_j) = \fr{\partial f(r,\zeta_j)}{\partial \theta} (\theta - \theta_j) \no (2.10)
$$
for some $\zeta_j \in [\theta_j,\theta_{j+1}], \ph j=\overline{0,N-1},$ provided $f(r,\theta)$ has no singularities on $C_r$
$$
f(r,\theta) - f(r_j,\theta) = \fr{\partial f(\nu_j,\theta)}{\partial r} (r-r_j) \no (2.11)
$$
for some $\nu_j \in [r_j,r_{j+1}], \ph j=\overline{1,N-1}.$ Assume $w_j\equiv r_je^{i\theta_j} \in KS(a,b,\la,\mu), \ph j=1,2.$ Then
$$
|w_1 - w_2| = \sqrt{(r_1\cos\theta_1 - r_2\cos\theta_2)^2 + (r_1\sin \theta_1 - r_2\sin \theta_2)^2}. \no (2.12)
$$
We have $|r_1\cos\theta_1 - r_22\cos\theta_2| \le r_1|\cos\theta_1 - \cos\theta_2| + |r_1 -r_2| |\cos\theta_j| \le r_1|\theta_1 -\theta_2| + |r_1 - r_2|$, and by analogy $|r_1\sin\theta_1 - r_2\sin\theta_2| \le  r_1|\theta_1 -\theta_2| + |r_1 - r_2|$. This together with (2.12) implies
$$
|w_1 - w_2| \le \sqrt{2}(r_1|\theta_1 - \theta_2| + |r_1 - r_2|).
$$
Therefore,
$$
KS(a,b,\la,\mu)\subset D(\La,\fr{\sqrt{2}}{2}(b(\mu-\la)+ b -a))  \no (2.13)
$$
for some $\La \in \mathbb{C}$.

We set $g(r,\theta)\equiv P(r e^{i\theta}), \ph h(r,\theta)\equiv P'(r e^{i\theta}). \ph r \ge 0, \ph 0 \le \theta \le 2\pi$.

\vskip10pt

{\bf Lemma 2.5.} {\it Let $x_s = \fr{2\pi s}{N}, \ph s =\overline{0,N}$. Then the following assertions are valid.

\noindent
(1) If $g(r,\theta_*) = 0$ for some $\theta_* \in [0,2\pi]$, then
$$
\min\limits_{s=\overline{0,N}}|g(r,x_s)| \le \Bigl(\sum\limits_{k=0}^{n-1}(n-k)|a_k| r^{n-k}\Bigr)\fr{2\pi}{N}. \no (2.14)
$$

\noindent
(2) If $h(r,\theta_*) = 0$ for some $\theta_* \in [0,2\pi]$, then
$$
\min\limits_{s=\overline{0,N}}|h(r,x_s)| \le \Bigl(\sum\limits_{k=0}^{n-2}(n-k)(n-k-1)|a_k| r^{n-k-1}\Bigr)\fr{2\pi}{N}.
$$
}

Proof. Let us prove (1). Assume $\theta_* \in [x_{s_0},x_{s_0+1}].$ Then by the Lagrange's formula $g(r,x_{s_0}) = \fr{\partial g(r,\widetilde{\theta})}{\partial \theta} (\theta_* - x_{s_0})$ for some $\widetilde{\theta} \in [x_{s_0},\theta_*]$. Hence,
$$
|g(r,x_{s_0}) \le \left|\fr{\partial g(r,\widetilde{\theta})}{\partial \theta}\right||\theta_* - x_{s_0}| \le  \left|\fr{\partial g(r,\widetilde{\theta})}{\partial \theta}\right|\fr{2\pi}{N}. \no (2.15)
$$
It is not difficult to verify that
$$
\left|\fr{\partial g(r,\widetilde{\theta})}{\partial \theta}\right| \le \sum\limits_{k=0}^{n-1}(n-k)|a_k| r^{n-k}.
$$
This together with (2.15) implies (2.14). The assertion (1) is proved. The assertion (2) can be proved by analogy with the proof of the assertion (1). The lemma is proved.

By analogy with the proof of Lemma 2.5 can be proved the following lemma

\vskip 10pt

{\bf Lemma 2.6.} {\it Let $r_s = a + \fr{(b-a) s}{N}, \ph s =\overline{0,N}$. Then the following assertions are valid.

\noindent
(3) If $h(r_*,\theta) = 0$ for some $r_*\ge 0$, then
$$
\min\limits_{s=\overline{0,N}}|h(r_s,\theta)| \le \Bigl(\sum\limits_{k=0}^{n-1}(n-k)|a_k| b^{n-k}\Bigr)\fr{b-a}{N}.
$$

\noindent
(4) If $h(r_*,\theta) = 0$ for some $r_*\ge 0$, then
$$
\min\limits_{s=\overline{0,N}}|h(r_s,\theta)| \le \Bigl(\sum\limits_{k=0}^{n-2}(n-k)(n-k-1)|a_k| b^{n-k-1}\Bigr)\fr{b-a}{N}.
$$
}

\phantom{aaaaaaaaaaaaaaaaaaaaaaaaaaaaaaaaaaaaaaaaaaaaaaaaaaaaaaaaaaaaaaaaa} $\Box$

\vskip 10pt

{\bf 3. \hskip 3pt An algorithm for finding the index of a function by a determined number of arithmetic operations.} Denote by $\Omega_{a,b}$ the set of all complex-valued continuous functions $F(x)$ on $[a,b]$, satisfying the conditions (see [6])

\noindent
(I) $F(x) \ne 0, \ph x \in [a,b]$,

\noindent
(II) $F(a) = F(b)$

\noindent
It is well known that for every $F \in\Omega_{a,b}$  the index $ind F(x)$ is well defined by the formula
$$
ind  \hskip 2pt F(x) \equiv \fr{1}{2\pi}[\arg F(x)]_{x=a}^b,
$$
where $[\arg F(x)]_{x=a}^b$ is the variation of $\arg F(x)$ for $x$ varying from $a$ to $b$. We set
$$
F_1(x)\equiv \fr{F(x)}{|F(x)|}\fr{F(a)}{|F(a)|}, \phh F\in \Omega_{a,b}.
$$
Let $a=x_1<x_2<\dots<x_N=b$ be a partition of the interval $[a,b]$ such that
$$
|F_1(x') - F_1(x'')|\le 2, \phh x',x'' \in [x_k,x_{k+1}], \phh k=\overline{1,N-1}. \no (3.1)
$$
and let $F_2(t)$ be  a function on $[a,b]$, which is a constant on each $[x_k,x_{k+1}] \ph (k=\overline{1,N-1})$  such that $F_2(x_k) = F_1(x_k), \ph k=\overline{1,N-1}$. Consider the systems of equations
$$
\sist{Re \hskip 2pt F_2(x) = 0,}{x\in\Delta_k}, \phh \sist{Im \hskip 2pt F_2(x) = 0,}{x\in\Delta_k}, \no (3.2)
$$
$k\in\{s \ph : \ph F_1(x_s) \ne F_1(x_{s+1}), \ph s=\overline{1,N-1},$ where $\Delta_k\equiv[x_k,x_{k=1})$ for $k=\overline{1,N-1}\}$ and $\Delta_N \equiv [x_{N-1},x_N].$ It was shown in [6] that these systems have finite number of solutions, called basic points of the function
$$
G(x) \equiv \fr{F_2(x)}{|F_2(x)|}.
$$
Let $a=\xi_1 <\dots < \xi_m=b$ be the basic points of $G(x)$. Consider the vector, called the indicator of $G(x)$(see [6]):
$$
I_G\equiv (G(\xi_1),\dots,G(\xi_m)) \phh (\mbox{obviously} \ph G(\xi_k) \in \{\pm1, \pm i\}, \ph k=\overline{1,m}).
$$
The number $m$ is called the length of $I_G$ and is denoted by $|I_G|$. If $|I_G|< 5$, then $ind F(x) = 0$ (see [6])
For calculating $ind F(x)$ for  $|I_G|\ge 5$ we use the following two types of reductions:

\noindent
1) If $G(\xi_p) = G(\xi_{p+1})$ for some $p$ then we rename $I_G \equiv (G(\xi_1),\dots,G(\xi_p),G(\xi_{p+2}),\dots,G(\xi_N))$, for which $|I_G| = N-1$.

\noindent
2) If $G(\xi_p) = G(\xi_{p+2})$ for some $p$ then we rename $I_G \equiv (G(\xi_1),\dots,G(\xi_p),G(\xi_{p+3}),\dots,G(\xi_N))$, for which $|I_G| = N-2$.

Using these reductions after finite number of reductions we obtain a vector of two types
$$
I_+\equiv (1, i, -1, -i, 1, i, \dots1, ), \phh I_- \equiv (1, -i, -1, i, 1, \dots,1)
$$
Then (see [6])
$$
ind \hskip 2pt F(x) = \pm \fr{|I_\pm| -1}{4},
$$
where $|I_\pm|$ is the length of $I_\pm$. It was shown in [6], that if $F(x)$ satisfies the Helder's condition
$$
|F(x') - F(x'')| \le M|x' - x''|^\alpha \phh (0 < \alpha \le 1) \no(3.3)
$$
and  $x_k = a + \fr{k-1}{N}(b-a), \ph k=\overline{1,N-1}$ for $N > (b-a)\left(\fr{2M}{m}\right)^{1/\alpha},$ where $m\equiv \min\limits_{k=\overline{1,N-1}}\{|F(x_k)|\}$, then the condition (3.1) holds. Thus for functions $F(x)$, satisfying the conditions (I), (II) and (3.3) we have a finite numerical method for calculation of $ind \hskip 2pt F(x)$. It was shown in [6] that
$$
I_G = \left(\fr{F_2(\xi_1)}{|F(\xi_1)|},\dots, \fr{F_2(\xi_q)}{|F(\xi_q)|}\right),
$$
where $\xi_1 < \dots < \xi_q$ are solutions of the systems
$$
\sist{\fr{Re \hskip 2pt F_1(x_{k+1}) - Re \hskip 2pt F_1(x_k)}{x_{k+1} - x_k} (x -x_k) = - Re \hskip 2pt F_1(x_k),}
{x\in\Delta_k},
$$
$$
\sist{\fr{Im \hskip 2pt F_1(x_{k+1}) - Im \hskip 2pt F_1(x_k)}{x_{k+1} - x_k} (x -x_k) = - Im \hskip 2pt F_1(x_k),}
{x\in\Delta_k},
$$
for all $k=\overline{1,N-1}$ for which $F_1(x_k) \ne F_1(x_{k+1})$.

{\bf 4. An algorithms for finding $ind \ph \fr{P'(z)}{P(z)}$ by a determined number of arithmetic operations}. For finding $\stackrel{ind}{_{C_r}}\fr{P'(z)}{P(z)}$ we use the following parametrization of $C_r: \ph z = r e^{i x}, \ph x \in[0,2\pi]$. Then
$$
\stackrel{ind}{_{C_r}}\fr{P'(z)}{P(z)}= \stackrel{ind}{_{x\in[0,2\pi]}}\fr{P'(r e^{ix})}{P(r e^{i x})}.
$$
If $r = a +  k \fr{b-a}{3}, \ph k=1,2$ where $K(a,b)$ is an annulus free of roots of $P(z)$and $P'(z)$, them by Lemma 2.1 $m\equiv \min\limits_{x\in[0,2\pi]}\left|\fr{P'(r e^{ix})}{P(r e^{i x})}\right|\ge  \frac{n}{(2\mathcal{R}_0)^n}\Bigl(\fr{b-a}{3}\Bigr)^{n-1}, \phh k=1,2.$ By Lemma 2.2 we have $M\equiv \max\limits_{x \in [0,2\pi]} \left|\fr{P'(r e^{ix})}{P(r e^{i x})}\right| \le \frac{9 n b}{(b -a)^2}, \ph r  = a + k\fr{b-a}{3}, \ph k=1,2,$ Then by (3.3) if we take $x_m = \fr{2\pi m}{N}, \ph m=0,1,\dots, N$ where
$N = \left[\fr{1}{2}\left(\fr{3}{\pi}\right)^{n+1}\mathcal{R}_0^{n+1}\right] + 1  \ge  \fr{4\pi M}{m}$ (here $[x]$ denotes the integer part of $x$) and use the algorithm of section 3 we can find $\stackrel{ind}{_{C_r}}\fr{P'(z)}{P(z)}$ by a determined number of arithmetic operations, provided $r = a +  k \fr{b-a}{3}, \ph k=1,2$ where $K(a,b)$ is an annulus free of roots of $P(z)$ and $P'(z)$. Let $\Gamma$ be the border of a segment of an annulus: $\Gamma\equiv \partial KS(a,b,\la,\mu)$. Then for finding $\stackrel{ind}{_{\small \Gamma}}\fr{P'(z)}{P(z)}$ we use the following parametrization of $\Gamma$
$$
\alpha(x)\equiv \left\{\begin{array}{l}
a e^{i(\la +\mu -x)}, \phh x \in[\la,\mu],\\
(x+a-\mu)e^{i\la}, \phh x\in [\mu,b-a+\mu],\\
b e^{i(x-b+a-\mu+\la)}, \phh[b-a+\mu,b-a+2\mu-\la],\\
(2 b - a +2\mu -\la - x) e^{i\mu}, \phh x \in [b-a+2\mu-\la,2 b - 2 a + 2\mu - \la].
 \end{array}\right.
$$
and take $\stackrel{ind}{_{\small \Gamma}}\fr{P'(z)}{P(z)}= \stackrel{ind}{_{x\in [\la,\nu]}} F(x),$ where $\nu\equiv 2 b - 2 a +\mu -\la, \ph F(x) \equiv \fr{P'(\al(x))}{P(\al(x))}$. It is not difficult to verify that
$$
F(x) = \left\{\begin{array}{l}
f(a,\la+\mu-x), \phh x \in[\la,\mu],\\
f(x+a-\mu,\la), \phh x\in [\mu,b-a+\mu],\\
f(b,x-b+a -\la +\mu), \phh [b-a+\mu,b-a+2\mu-\la],\\
f(2 b - a +2\mu -\la - x,\mu), \phh  x \in [b-a+2\mu-\la,2 b - 2 a + 2\mu - \la]
\end{array}
\right.
$$
Let $KS(a-c,b+c,\la-\theta,\mu+\theta)$ be an annulus free of roots of $P(z)$ and $P'(Z)$. Then using Lemmas 2.1 - 2.4 one can show that $m\equiv \min\limits_{x\in\Gamma}F(x) \ge \fr{n}{(2 \mathcal{R}_0)^n}\min\{\left(\fr{b-a}{3}\right)^n,\left(\min\{c,a \tan \theta\}\right)^{n-1}\}$ and $m\equiv \max\limits_{x\in\Gamma}F(x) \le n \max\left\{\fr{9b}{(b-)^2}, \fr{1}{\left(\min\{c,a \tan \theta\}\right)^2}\right\}$.
Then by (3.3) if we take $x_m = \fr{2\pi m}{N}, \ph m=0,1,\dots, N$, where
$$
N =\left[\fr{2(2\mathcal{R}_0)^n \max\left\{\fr{9b}{(b-a)^2}, \fr{1}{\left(\min\{c,a \tan \theta\}\right)^2}\right\}}{\min\{\left(\fr{b-a}{3}\right)^n,\left(\min\{c,a \tan \theta\}\right)^{n-1}\}}\right] + 1 \ge \fr{2\nu M}{m},
$$
and use the algorithm of section 3 we can find $\stackrel{ind}{_{x\in\Gamma}}F(x)$ by a determined number of arithmetic operations.

\vskip 10pt

{\bf 5.  The separation   of an annulus with  determined safe (neutral) zones by a determined number of arithmetic operations.}

{\bf Definition 5.1} {\it The annuluses $K(a,a+c_1)$ and $K(b-c_2,b) \ph (a+c_1 \le b -c_2)$ are called safe (neutral) zones for an annulus $K(a,b)$, if $K(a,b) \backslash  K(a+c_1,b - c_2)$ is free of roots of $P(z)$ and $P'(z)$.
}
 By results  of section 4 if $K(a,b)$ is an annulus with determined safe (neutral) zones $K(a,a+c_1)$ and $K(b-c_2,b)$ ($c_1$ and $c_2$ are determined [established]), then the quantity $n_{a,b}$ of roots of $P(z)$,  contained  in  $K(a,b)$,  can be found by determined number of arithmetic operations by using
the formula
$$
n_{a,b} = \stackrel{ind}{_{C_{r_2}}}\fr{P'(z)}{P(z)} - \stackrel{ind}{_{C_{r_1}}}\fr{P'(z)}{P(z)}, \ph r_1 \equiv a +c_1/2, \ph r_2 \equiv b - c_2/2. \no (5.1)
$$
In this section we show how one can separate an annulus with the determined safe zones from the another annulus with the determined safe zones. Let $K(a,b)$ be an annulus with the determined safe zones. Let us break it into $2n$ annuluses with equal widths: $K(a,b) = \cup_{m=0}^{2n-1}K(r_m,r_{m+1})$, where $r_m\equiv a + \fr{m(b-a)}{2n}, \ph m=\overline{0,2n}$. Denote by $n'_{a,b}$ the quantity of roots of $P'(z)$ counted as many times as their multiplicity.
Since, obviously, $n_{a,b} + n'_{a,b} \le 2n-1$, at least one of $\stackrel{_\circ}{K}(r_m,r_{m+1}), \ph m=\overline{0,2n-1}$ is free of roots of $P(z)$ and $P'(z)$. Let $\stackrel{_\circ}{K}(r_{m_0},r_{m_0+1})$ be one of such annuluses. Then
$$
\min\limits_{z\in C_{\et_{m_0}}}|P(z)| \ge |a_0|\Bigl(\fr{b-a}{6n}\Bigr)^n, \no (5.2)
$$
$$
\min\limits_{z\in C_{\zeta_{m_0}}}|P(z)| \ge |a_0|\Bigl(\fr{b-a}{6n}\Bigr)^n, \no (5.3)
$$

$$
\min\limits_{z\in C_{\et_{m_0}}}|P'(z)| \ge n|a_0|\Bigl(\fr{b-a}{6n}\Bigr)^{n-1}, \no (5.4)
$$
$$
\min\limits_{z\in C_{\zeta_{m_0}}}|P'(z)| \ge n|a_0|\Bigl(\fr{b-a}{6n}\Bigr)^{n-1}, \no (5.5)
$$
where $\et_{m}\equiv r_{m} + \fr{b-a}{6n}, \ph \zeta_{m}\equiv r_{m} + \fr{b-a}{3n}, \ph m=\overline{0,2n-1}$. Let $x_s = \fr{2\pi s}{N}, \ph s=\overline{0,N}$. If $P(z_*) = 0$ for some $z_*\in C_{\et_m}$ (for some $z_*\in C_{\zeta_m}$), then by Lemma 2.5  (1)
$$
\min\limits_{s=\overline{0,N}}|P(\et_m e^{i x_s})| \le  \Bigl(\sum\limits_{k=0}^{n-1}(n-k)|a_k| \et_m^{n-k}\Bigr)\fr{2\pi}{N}, \no (5.6)
$$
$$
\left(\min\limits_{s=\overline{0,N}}|P(\et_m e^{i x_s})| \le  \Bigl(\sum\limits_{k=0}^{n-1}(n-k)|a_k| \zeta_m^{n-k}\Bigr)\fr{2\pi}{N}\right), \no (5.7)
$$
and if  $P'(z_*) = 0$ for some $z_*\in C_{\et_m}$ (for some $z_*\in C_{\zeta_m}$), then by Lemma 2.5  (2)
$$
\min\limits_{s=\overline{0,N}}|P'(\et_m e^{i x_s})| \le  \Bigl(\sum\limits_{k=0}^{n-2}(n-k)(n-k-1)|a_k| \et_m^{n-k-1}\Bigr)\fr{2\pi}{N}, \no (5.8)
$$
$$
\left(\min\limits_{s=\overline{0,N}}|P'(\et_m e^{i x_s})| \le   \Bigl(\sum\limits_{k=0}^{n-2}(n-k)(n-k-1)|a_k| \zeta_m^{n-k-1}\Bigr)\fr{2\pi}{N}\right). \no (5.9)
$$
It is not difficult to verify that if we take
$$
N=\left [\fr{2\pi(6n)^{n-1}}{|a_0|(b-a)^{n-1}}\\max\left\{\fr{\sum\limits_{k=0}^{n-1}(n-k)|a_k| b^{n-k}}{(b-a)}, \fr{\sum\limits_{k=0}^{n-2}(n-k)(n-k-1)|a_k| b^{n-k-1}}{n}\right\}\right] +1, \no (5.10)
$$
then from (5.2)--(5.9) we obtain that for any $m$ at least one of  the circles $C_{\et_m}$ and $C_{\zeta_m}$ is not free of roots of $P(z)$ and $P'(z)$ if at least one of the following inequalities hold.
$$
\min\limits_{s=\overline{0,N}}P(\et_m e^{ix_s})| < |a_0|\Bigl(\fr{b-a}{6n}\Bigr)^n, \no (5.11)
$$
$$
\min\limits_{s=\overline{0,N}}P(\zeta_m e^{ix_s})| < |a_0|\Bigl(\fr{b-a}{6n}\Bigr)^n, \no (5.12)
$$
$$
\min\limits_{s=\overline{0,N}}P'(\et_m e^{ix_s})| < |a_0|\Bigl(\fr{b-a}{6n}\Bigr)^{n-1}, \no (5.13)
$$
$$
\min\limits_{s=\overline{0,N}}P'(\zeta_m e^{ix_s})| < |a_0|\Bigl(\fr{b-a}{6n}\Bigr)^{n-1}. \no (5.14)
$$
Now we have all data for describing an algorithm for separation (selection) of an annulus $K(a_1,b_1)$ with the safe zones from a given annulus with safe zones by determined number of arithmetic operations. The algorithm is the following

\noindent
1. step 1. take  $\et_{m}\equiv r_{m} + \fr{b-a}{6n}, \ph \zeta_{m}\equiv r_{m} + \fr{b-a}{3n}$, where $r_m\equiv a + \fr{m(b-a)}{2n}, \ph m=\overline{0,2n-1}$,

\noindent
2. step 2. calculate (find) $N$ by the formula (5.10),

\noindent
3. step 2. chose $m=m_0$ such that all inequalities (5.11)--(5.14) are not satisfied (by (5.2)--(5.5) at least one such $m_0$ exists),

\noindent
4. step 4. calculate $n_{a,\et_{m_0}}$ and $a_{\zeta_{m_0},b}$ by the formula (5.1)

\noindent
5. if $n_{a,\et_{m_0}} = 0$, then we put $K(a_1,b_1) \equiv K(\fr{\et_{m_0} + \zeta_{m_0}}{2},b)$ (which have safe zones);
if $n_{a,\et_{m_0}} \ne0$, then we put $K(a_1,b_1) \equiv K(a,\fr{\et_{m_0} + \zeta_{m_0}}{2})$ (which have also safe zones).
 \vskip 10pt

{\bf 6. The separation (selection) of a segment with a determined safe (neutral) zone  from an annulus with determined safe zones.}

{\bf Definition 6.1.} {\it A segment $KS(a,b,\la,\mu)$ is called a segment with a safe (neutral) zone, if there exist $c_1> 0, c_2> 0, \theta_1> 0, \theta_2> 0$ such that $a+c_1 \le b-c_2, \ph \la + \theta_1 \le \mu -\theta_2$ and
$KS(a,b,\la,\mu) \backslash KS(a+c_1,b-c_2,\la+\theta_1,\mu-\theta_2)$ is free of roots of $P(z)$ and $P'(z)$. Under these conditions the set $KS(a,b,\la,\mu) \backslash KS(a+c_1,b-c_2,\la+\theta_1,\mu-\theta_2)$ is called a safe (neutral) zone for $K(a,b,\la,\mu)$.
}

Let $K(a,b)$ be an annulus with determined safe zones $K(a,a+c_1)$ and $K(b-c_2,b)$. We break $K(a,b)$ into $2n$ congruent segments $KS(a,b,\theta_m,\theta_{m+1})$, where $\theta_m \equiv \fr{\pi m}{n}, \ph m=\overline{0,2n}.$ Then there exists $m=m_0$ such that $\stackrel{_\circ}{K}S(a,b,\theta_{m_0},\theta_{m_0+1})$ is  free of roots of $P(z)$ and $P'(z)$. Let $\et_m=\theta_m + \fr{\pi}{3n}, \zeta_m = \theta_m + \fr{2\pi}{3n}, \ph m=\overline{0,2n-1}$. Then it is not difficult to  verify that
$$
\min\limits_{a\le r \le b}|P(r e^{i \et_m})| \ge |a_0|\Bigl(\sin\fr{\pi}{3n}\Bigr)^n, \no (6.1)
$$
$$
\min\limits_{a\le r \le b}|P(r e^{i \zeta_m})| \ge |a_0|\Bigl(\sin\fr{\pi}{3n}\Bigr)^n, \no (6.2)
$$
$$
\min\limits_{a\le r \le b}|P'(r e^{i \et_m})| \ge n|a_0|\Bigl(\sin\fr{\pi}{3n}\Bigr)^{n-1}, \no (6.3)
$$
$$
\min\limits_{a\le r \le b}|P'(r e^{i \zeta_m})| \ge n|a_0|\Bigl(\sin\fr{\pi}{3n}\Bigr)^{n-1}. \no (6.4)
$$
We set $r_z=a + \fr{(b-a)s}{N}, \ph s =\overline{0,N}.$ If $P(r_* e^{i\et_m})=0 \ph (P(r_* e^{i\zeta_m})=0)$ for some $r_*\in[a,b]$ and for some $m$, then by Lemma 2.6 (3)
$$
\min\limits_{s=\overline{0,N}}|P(r_s e^{i\et_m})| \le \Bigl(\sum\limits_{k=0}^{n-1}(n-k)|a_k| b^{n-k}\Bigr)\fr{b-a}{N}, \no (6.5)
$$
$$
\left(\min\limits_{s=\overline{0,N}}|P(r_s e^{i\zeta_m})| \le \Bigl(\sum\limits_{k=0}^{n-1}(n-k)|a_k| b^{n-k}\Bigr)\fr{b-a}{N}\right), \no (6.6)
$$
and if $P'(r_* e^{i\et_m})=0 \ph (P'(r_* e^{i\zeta_m})=0)$ for some $r_*\in[a,b]$ and for some $m$, then by  Lemma~ 2.6 (4)
$$
\min\limits_{s=\overline{0,N}}|P'(r_s e^{i\et_m})| \le \Bigl(\sum\limits_{k=0}^{n-2}(n-k)(n-k-1)|a_k| b^{n-k-1}\Bigr)\fr{b-a}{N}, \no (6.7)
$$
$$
\left(\min\limits_{s=\overline{0,N}}|P'(r_s e^{i\zeta_m})| \le \Bigl(\sum\limits_{k=0}^{n-2}(n-k)(n-k-1)|a_k| b^{n-k-1}\Bigr)\fr{b-a}{N}\right), \no (6.8)
$$
On the basis (6.1)--(6.8)  one can easily  verify that if we take
$$
N=\left[\fr{(b-a)b^{n-k-1}}{|a_0|\Bigl(\sin\fr{\pi}{3n}\Bigr)^{n-1}}\max\left\{\fr{\Bigl(\sum\limits_{k=0}^{n-1}(n-k)|a_k| b\Bigr)}{\sin\fr{\pi}{3n}}, \fr{\sum\limits_{k=0}^{n-2}(n-k)(n-k-1)|a_k|}{n}\right\}\right] + 1, \no (6.9)
$$
then the lines $\theta = \et_m, \ph \theta = \zeta_m, \ph a \le r \le b$ are free of roots of $P(z)$ and $P'(z)$ provided
$$
\min\limits_{s=\overline{0,N}}|P(r_s e^{i\et_m})| \ge |a_0|\Bigl(\sin\fr{\pi}{3n}\Bigr)^n, \no (6.10)
$$
$$
\min\limits_{s=\overline{0,N}}|P(r_s e^{i\zeta_m})| \ge |a_0|\Bigl(\sin\fr{\pi}{3n}\Bigr)^n, \no (6.11)
$$
$$
\min\limits_{s=\overline{0,N}}|P'(r_s e^{i\et_m})| \ge n |a_0|\Bigl(\sin\fr{\pi}{3n}\Bigr)^{n-1}, \no (6.12)
$$
$$
\min\limits_{s=\overline{0,N}}|P'(r_s e^{i\zeta_m})| \ge n |a_0|\Bigl(\sin\fr{\pi}{3n}\Bigr)^{n-1}. \no (6.13)
$$
Now to separate a segment $KS(a,b,\la,\mu)$  with a determined safe zone from an annulus $K(a,b)$ with determined safe zones by determined number of arithmetic operations we can follow the following steps of actions.

\noindent
1. Step 1. Find $N$ by formula (6.9)

\noindent
2. Step 2. Set $et_m=\theta_m +\fr{\pi}{3n}, \ph \zeta_m= \theta_m + \fr{2\pi}{3n},$ where $\theta_m = \fr{\pi m}{n}, \ph m=\overline{0,2n-1}.$

\noindent
3. Step 3. Find $m=m_0$ for which the inequalities (6.10)--(6.13) are satisfied (due to (6.1)--(6.4) at least one of such $m_0$ exists)

\noindent
4. Step 4. Set $KS(a,b,\la,\mu \equiv KS(a,b,\zeta_{m_0},2\pi+\et_{m_0})$ (this is a segment with a determined safe zone).

\noindent
5. Step 5. Calculate $\stackrel{ind}{_{\Gamma_{m_0}}}\fr{P'(z)}{P(z)}$ according to the method of section 4, where $\Gamma_{m_0} \equiv \partial KS(a,b,\zeta_{m_0},2\pi+\et_{m_0})$.

\noindent
6. Step 6. Calculate $n_{a,b}$ according to the method of section 4.

\noindent
If $n_{a,b} \ne \stackrel{ind}{_{\Gamma_{m_0}}}\fr{P'(z)}{P(z)}$, then go to step 3 and find another $m_0$.

{\bf 7. The separation (selection) of a segment with a determined safe (neutral) zone  from another segment  with a determined safe zone.} Let $KS(a,b,\la,\mu)$ be a segment with a determined safe zone.

We break $K(a,b,\la,\mu)$ into $2n$ congruent segments $KS(a,b,\theta_m,\theta_{m+1})$, where $\theta_m \equiv \la + \fr{(\mu -\la) m}{2n}, \ph m=\overline{0,2n}.$ Then there exists $m=m_0$ such that $\stackrel{_\circ}{K}S(a,b,\theta_{m_0},\theta_{m_0+1})$ is  free of roots of $P(z)$ and $P'(z)$. Let $\et_m=\theta_m + \fr{\mu-\la}{6n}, \ph  \zeta_m = \theta_m + \fr{\mu-\la}{3n}, \ph m=\overline{0,2n-1}$. Then it is not verify that
$$
\min\limits_{a\le r \le b}|P(r e^{i \et_m})| \ge |a_0|\Bigl(\sin\fr{\mu-\la}{6n}\Bigr)^n, \no (6.14)
$$
$$
\min\limits_{a\le r \le b}|P(r e^{i \zeta_m})| \ge |a_0|\Bigl(\sin\fr {\mu-\la}{6n}\Bigr)^n, \no (6.15)
$$
$$
\min\limits_{a\le r \le b}|P'(r e^{i \et_m})| \ge n|a_0|\Bigl(\sin\fr{\mu-\la}{6n}\Bigr)^{n-1}, \no (6.16)
$$
$$
\min\limits_{a\le r \le b}|P'(r e^{i \zeta_m})| \ge n|a_0|\Bigl(\sin\fr{\mu-\la}{6n}\Bigr)^{n-1}. \no (6.17)
$$
We set $r_z=a + \fr{(b-a)s}{N}, \ph s =\overline{0,N}.$ If $P(r_* e^{i\et_m})=0 \ph (P(r_* e^{i\zeta_m})=0)$ for some $r_*\in[a,b]$ and for some $m$, then by Lemma 2.6 (3)
$$
\min\limits_{s=\overline{0,N}}|P(r_s e^{i\et_m})| \le \Bigl(\sum\limits_{k=0}^{n-1}(n-k)|a_k| b^{n-k}\Bigr)\fr{b-a}{N}, \no (6.18)
$$
$$
\left(\min\limits_{s=\overline{0,N}}|P(r_s e^{i\zeta_m})| \le \Bigl(\sum\limits_{k=0}^{n-1}(n-k)|a_k| b^{n-k}\Bigr)\fr{b-a}{N}\right), \no (6.19)
$$
and if $P'(r_* e^{i\et_m})=0 \ph (P'(r_* e^{i\zeta_m})=0)$ for some $r_*\in[a,b]$ and for some $m$, then by  Lemma~ 2.6 (4)
$$
\min\limits_{s=\overline{0,N}}|P'(r_s e^{i\et_m})| \le \Bigl(\sum\limits_{k=0}^{n-2}(n-k)(n-k-1)|a_k| b^{n-k-1}\Bigr)\fr{b-a}{N}, \no (6.20)
$$
$$
\left(\min\limits_{s=\overline{0,N}}|P'(r_s e^{i\zeta_m})| \le \Bigl(\sum\limits_{k=0}^{n-2}(n-k)(n-k-1)|a_k| b^{n-k-1}\Bigr)\fr{b-a}{N}\right). \no (6.21)
$$
Using (5.14)--(6.21) one can easily show that  if we take
$$
N=\left[\fr{(b-a)b^{n-k-1}}{|a_0|\Bigl(\sin\fr{\mu-\la}{6n}\Bigr)^{n-1}}\max\left\{\fr{\Bigl(\sum\limits_{k=0}^{n-1}(n-k)|a_k| b\Bigr)}{\sin\fr{\mu-\la}{6n}}, \fr{\sum\limits_{k=0}^{n-2}(n-k)(n-k-1)|a_k|}{n}\right\}\right] + 1, \no (6.22)
$$
then the lines $\theta = \et_m, \ph \theta = \zeta_m, \ph a \le r \le b$ are free of roots of $P(z)$ and $P'(z)$ provided
$$
\min\limits_{s=\overline{0,N}}|P(r_s e^{i\et_m})| \ge |a_0|\Bigl(\sin\fr{\mu-\la}{6n}\Bigr)^n, \no (6.23)
$$
$$
\min\limits_{s=\overline{0,N}}|P(r_s e^{i\zeta_m})| \ge |a_0|\Bigl(\sin\fr{\mu-\la}{6n}\Bigr)^n, \no (6.24)
$$
$$
\min\limits_{s=\overline{0,N}}|P'(r_s e^{i\et_m})| \ge n |a_0|\Bigl(\sin\fr{\mu-\la}{6n}\Bigr)^{n-1}, \no (6.25)
$$
$$
\min\limits_{s=\overline{0,N}}|P'(r_s e^{i\zeta_m})| \ge n |a_0|\Bigl(\sin\fr{\mu-\la}{6n}\Bigr)^{n-1}. \no (6.26)
$$
Now to separate a segment $K(a,b,\la_1,\mu_1)$   with a determined safe zone from another segment $K(a,b,\la,\mu)$ with a determined safe zone by determined number of arithmetic operations we can follow the following steps of actions.

\noindent
1. Step 1. Find $N$ by formula (6.22)

\noindent
2. Step 2. Set $\et_m=\theta_m +\fr{\mu-\la}{6n}, \ph \zeta_m= \theta_m + \fr{\mu-\la}{3n},$ where $\theta_m = \fr{\pi m}{n}, \ph m=\overline{0,2n-1}.$

\noindent
3. Step 3. Find $m=m_0$ for which the inequalities (6.23)--(6.26) are satisfied (due to (6.14)--(6.17) at least one of such $m_0$ exists)

\noindent
4. Step 4. Set $KS(a_1,b_1,\la_1,\mu_1) \equiv KS(a,b,\zeta_{m_0},2\pi+\et_{m_0})$ (this is a segment with a determined safe zone).

\noindent
5. Step 5. Calculate $\stackrel{ind}{_{\Gamma_{m_0}}}\fr{P'(z)}{P(z)}$ according to the method of section 4, where $\Gamma_{m_0} \equiv \partial KS(a,b,\zeta_{m_0},2\pi+\et_{m_0})$.

\noindent
6. Step 6. Calculate $n_{a,b}$ according to the method of section 4.

\noindent
If $n_{a,b} \ne \stackrel{ind}{_{\Gamma_{m_0}}}\fr{P'(z)}{P(z)}$, then go to step 3 and find another $m_0$.

\vskip 10pt

{\bf 8. The algorithm of localization of roots of $P(z)$ by determined number of arithmetic operations.} According to Theorem 2.1 and Corollary 2.1 all roots of $P(z)$ are located in the annulus $K(\rho_0,\mathcal{R}_0)$ and all roots of $P'(z)$ are located in the annulus $K(\rho_1,\mathcal{R}_1)$, where
$$
\rho_1\equiv \fr{|a_{n-1}|}{|a_{n-1}| + \max\{n|a_0|,\dots,2|a_{n-2}|\}}, \ph \mathcal{R}_1\equiv \fr{\max\{(n-1)|a_1|,\dots,|a_{n-1}|\}}{n|a_0|}\le \mathcal{R}_0,
$$
provided $a_{n-1} \ne 0$. If $a_{m+1} = \dots = a_{n-1} = 0$ and $a_m\ne 0$, then all roots of $P'(z)$ except $z_0=0$ are located in $K(\rho_2,\mathcal{R}_1)$, where
$$
\rho_2\equiv \fr{(n-m)|a_m|}{(n-m)|a_m| + \max\{n|a_0|,\dots,(n-m)|a_m|\}}.
$$
We set $\widetilde{\rho}\equiv \sist{\rho_1/2, \ph if \ph a_{n-1}\ne 0,}{\rho_2/2, \ph if \ph a_m\ne0, \ph a_{m+1}=\dots=a_{n_1}=0.}$ Then $K(\widetilde{\rho},\max\{\mathcal{R}_0,\mathcal{R}_1\}+1) = K(\widetilde{\rho},\mathcal{R}_0+1)$ is an annulus with determined safe zones. For localization of roots of $P(z)$ by finite number of arithmetic operations we can now use the following algorithm.

\noindent
1. Step 1. Chose $\varepsilon < \varepsilon/4$, where $\varepsilon_0=sep(P)$ is defined in section 2.

\noindent
2. Step .2. Using the algorithm of section 6 separate from $K(\widetilde{\rho},\mathcal{R}_0+1)$ all annuluses $K(a_m,b_m), \ph m=\overline{1,n_0}$ with $b_m-a_m \le \fr{\sqrt{2}\varepsilon}{\mathcal{R}_0+2}, \ph m=\overline{1,n_0}$ and with determined safe zones, containing all roots of $P(z)$ (by determined number of arithmetic operations).

\noindent
2. Step 2. Using the algorithms of sections 6 and 7 separate (by determined number of arithmetic operations) from all  annuluses $K(a_m,b_m), \ph m=\overline{1,n_0}$
all segments $KS(a_m,b_m,\la_{mk},\mu_{mk}), \ph k=\overline{1,l_m},\ph m=\overline{1,n_0}$ with $\mu_{mk} - \la_{mk} \le \fr{\sqrt{2}\varepsilon}{\mathcal{R}_0+2}, \ph k=\overline{1,l_m},\ph m=\overline{1,n_0}$, containing all roots of $P(z)$..

After realization  of this algorithm we obtain $n=l_1+\dots+l_{n_0}$ segments, each of which contain exactly one root of $P(z)$.  By (2.13) each of these segments contains in a disc of $\varepsilon$ radius.

Show that the mentioned above algorithm is realizable by determined number of arithmetic operations. Starting with $K(\widetilde{\rho},\mathcal{R}_0+1)$ using the algorithm of section 6 by determined number of arithmetic operations we can separate an annulus $K(a_1,b_1)$, containing some roots of $P(z)$ with $h_1\equiv b_1-a_1  \le \mathcal{R}_0(1- \fr{1}{2n})$. From $K(a_1,b_1)$ we can separate $K(a_2,b_2)$, containing some roots of $P(z)$ with $h_2\equiv b_1-a_1  \le h_1(1- \fr{1}{2n}) \le \mathcal{R}_0(1- \fr{1}{2n})^2$, so on. After N steps of such separations we obtain $K(a_N,b_N)$, containing some roots of $P(z)$, for which $b_N-a_N\le \mathcal{R}_0(1-\fr{1}{2n})^N$. Therefore,
$$
N\le \left[\fr{\ln (\sqrt{2} \hskip 3pt \varepsilon) - \ln(\mathcal{R}_0(\mathcal{R}_0 +2))}{\ln(2n-1) - \ln(2n)}\right] +1.
$$
Thus after determined number of arithmetic operations we can separate all annuluses with required width, containing all roots of $P(z)$. By analogy it can be sown that from the separated annuluses it can be separated discs with the required length containing all roots of $P(z)$ by definite number of arithmetic operations.

\vskip 10pt

\centerline{\bf References}

\vskip 10pt

\noindent
1. M. Mignote, Strasbourg, Some Useful Bounds. Computing, Suppl. 4, 259--263 (1982).

\noindent
2. A. G. Akristas, A. W. Strzebovitski, P. S. Vigklas, Improving the Performance of the \linebreak \ph  Continued Fractions Method Using New Bounds of Positive Roots. Nonlinear Analysis: \linebreak \ph  Modeling and Control, 2008, vol. 13, num. 3, pp. 265--279.

\noindent
3. Chu Keng Yap, Functional Problems in Algorithmic Algebra. Courant Institute of \linebreak \ph   Mathematical   Sciences, New York University, 251, Mercer Street, New York, \linebreak \ph   NY 10012, 1993.

\noindent
4. A. Edelman and E. Kostlan, How many Zeroes of a Random Polynomial Are real?  \linebreak \ph  Bulletin (New Series) of rhe AMS, vol. 30, num. 1, 1996, pp. 1--27.

\noindent
5. B. L. van der Varden, Algebra, Moskow, "Nauka", \ph 1979 (the russian translation of \linebreak \ph  B. L. VAN DER VARDEN ALGEBRA I, ACHTE AUFLING DER MODERNEN \linebreak \ph   ALGEBRA, Springer verlag, Berlin, Haidelberg, new York, 1971).

\noindent
6. G. A. Grogorian, On a numerical method for calculation of the index of a function. \linebreak \ph  Proceedings  of the NAS of Armenia, Mathematica, vol. 40. num. 1, 2005, pp. 47--60.

\end{document}